\documentclass{amsart}

\usepackage{amscd}
\usepackage{amssymb, amsfonts}
\usepackage{amsrefs}
\input xy
\xyoption{all}
\usepackage{color}
\usepackage[T1]{fontenc}
\usepackage{tikz}
\usepackage{mathrsfs}

\definecolor{bettergreen}{rgb}{0,.7,0}

\newcommand{\C}{{\mathbb C}}
\newcommand{\R}{{\mathbb R}}
\newcommand{\Z}{{\mathbb Z}}

\def\CC{{\mathbb C}}
\def\ZZ{{\mathbb Z}}

\newcommand{\cO}{\mathcal{O}}

\newcommand{\bF}{\mathbf{F}}
\newcommand{\ri}{\mathcal{O}}

\newcommand{\fg}{\mathfrak{g}}

\newcommand{\fsp}{\mathfrak{sp}}
\newcommand{\fsl}{\mathfrak{sl}}

\newcommand{\reg}{\mathrm{reg}}
\newcommand{\ur}{\mathrm{ur}}

\newcommand{\bG}{\mathbf{G}}

\newcommand{\oi}{{\bf \mathrm{O}}}

\newcommand{\rf}{k}

\newcommand\ldpo[1][\ri]{{\mathcal L}_{#1}}

\newcommand\cA{{\mathcal A}}
\newcommand\cB{{\mathcal B}}
\newcommand\ord{\mathrm{ord}}
\newcommand\ac{\overline{\mathrm{ac}}}

\newcommand\cC{{\mathcal C}}

\newcommand{\scD}{\mathscr{D}}

\newcommand{\K}{F}
\newcommand{\mexp}{\mathbf{e}}

\def\llp{\mathopen{(\!(}}

\def\rrp{\mathopen{)\!)}}

\theoremstyle{plain}
\newtheorem{thm}{Theorem}
\newtheorem{theorem}[thm]{Theorem}

\theoremstyle{definition}

\newtheorem{defn}[thm]{Definition}

\title[]{Motivic functions, integrability, and uniform in $p$ bounds for orbital integrals}
\author{Raf Cluckers, Julia Gordon and Immanuel Halupczok}

\begin{document}

\maketitle
\section{Introduction}
The purpose of this short article is two-fold: to announce the new results related to analytic properties of  motivic (exponential) functions, which we hope could be widely applicable; and to summarize the recent applications of these results in harmonic analysis on $p$-adic groups.

{Concretely}, motivic (exponential) functions are complex-valued functions on $p$-adic manifolds defined uniformly in $p$
by means of a rich yet sufficiently manageable language of logic. The word ``exponential'' refers to the situation when additive characters of the local field (leading to exponential sums) are included in the set of ingredients one is allowed to use when making such a function.
The class of motivic exponential functions has two key features that make it particularly  useful. On the one hand, these functions  can be described in a field-independent manner, and thus are an ideal tool  for transferring results between different local fields. On the other hand, the class of motivic (exponential) functions is closed under integration with respect to parameters, and thanks to this property one can show that many functions naturally arising in harmonic analysis on $p$-adic groups belong to this class.

The project of using motivic integration in harmonic analysis on $p$-adic groups was initiated by T.C. Hales in 1998, and {goes} under the emerging name of ``motivic harmonic analysis''. One of the goals of this article is to collect results about ``definability'' of various objects arising in representation theory of $p$-adic groups and the applications in harmonic analysis that follow from such definability,
summarizing the state {of the art} of this project.

In the following sections, we {summarize} most of the definitions we use that relate to (motivic) functions, motivic measures, etc. However, we do not {include} here any definitions related to harmonic analysis on reductive groups, as they are widely available in other sources, for example, in the beautiful article \cite{kottwitz:clay}.

{\bf Acknowledgment.} This article summarizes the results of a project carried out over several years. It certainly would not have happened without the influence of Tom Hales and the prior work of Jonathan Korman and Jyotsna Diwadkar.
It is a pleasure to thank many people who have helped us at various stages of the project, (and those who collaborated on it in earlier stages), in particular, Jan Denef, Fran\c cois Loeser,  Fiona Murnaghan, Clifton Cunningham, Loren Spice,
Sug-Woo Shin, Nicolas Templier,  Tasho Statev-Kaletha, William Casselman, Gopal Prasad, and Lance Robson. I.H. was supported by the SFB~878 of the Deutsche
Forschungsgemeinschaft; J.G. was supported by NSERC; R.C. was supported in part by the European Research Council under the European Community's Seventh Framework Programme (FP7/2007-2013) with ERC Grant Agreement nr. 246903 NMNAG, by the Labex CEMPI  (ANR-11-LABX-0007-01), and by the Fund for Scientific Research of Flanders, Belgium (grant G.0415.10).

\section{What is a motivic (exponential) function?}
Informally, motivic functions are built from definable functions in the Denef-Pas language and motivic \emph{exponential} functions are built from motivic functions, additive characters on local fields, and definable functions serving as arguments of the characters. As such, they are given independently of the local field and can be interpreted in any local field and  with any additive character. 
Let us recall the definition of the Denef-Pas language first, and proceed to motivic (exponential) functions afterwards, where we {slightly} simplify the terminology and notation of \cite{CGH-1}.

\subsection{Denef-Pas language}\label{sub:DP}
The Denef-Pas language is a first order language of logic designed for working
with valued fields. Formulas in this language will allow us to uniformly {handle} sets and functions for all {local} fields.
We start by defining two sublanguages of the language of
Denef-Pas: the language of rings and the Presburger language.
\subsubsection{The language of rings}

Apart from the symbols for variables $x_1, \dots, x_n,\dots$ and the
usual logical symbols equality `$=$', parentheses `$($', `$)$', the quantifiers `$\exists$', `$\forall$', and the logical operations conjunction `$\wedge$', negation `$\neg$', disjunction `$\vee$', the language of rings consists of the following symbols:
\begin{itemize}
\item constants `$0$', `$1$';
\item binary functions `$\times$', `$+$'.
\end{itemize}

A (first-order) formula in the language of rings is any syntactically correct formula
built out of these symbols. One usually omits the words `first order'.

If a formula in the language of rings has $n$ free (i.e.\ unquantified) variables then it defines a subset of $R^n$ for any ring $R$.
For example the formula ``$\exists x_2\, (x_2\times x_1 = 1)$'' defines the set of units $R^\times$ in any ring $R$; note that
quantifiers {(by convention)} always run over the ring in question.
Note also that quantifier-free formulas in the language of rings define constructible sets (in the sense of algebraic geometry).

\subsubsection{Presburger language}\label{Pres}
A formula in Presburger's language is built out of variables running over $\Z$, the logical symbols (as above) and
symbols `$+$', `$\le$', `$0$', `$1$', and for each $d=2,3,4,\dots$, a symbol `$\equiv_d$' to denote the binary
relation $x\equiv y \pmod{d}$.
Note the absence of the symbol for multiplication.

Since multiplication is not allowed, sets and functions defined by formulas in the Presburger language are in fact very basic, cf.~\cite{CPres} or \cite{Presburger}. For example, definable subsets of the line $\ZZ$ are finite unions of arithmetic progressions (in positive or negative direction) and points.
{By \cite{Presburger}, quantifiers are never needed to describe Presburger sets, so}
all definable sets are of simple form.

\subsubsection{Denef-Pas language}\label{DP}
The Denef-Pas language is a three-sorted language in the sense that its formulas speak about
three different ``sorts'' of elements: those of the valued field, of the residue field and of the value group (which will always be $\Z$ in our setting).
Each variable in such a formula only runs over the elements of one of the sorts,
i.e., there are three disjoint sets of symbols for the variables of the different sorts.
For a formula to be syntactically correct, one has to pay attention to the sorts when composing functions
and plugging them into relations.

In addition to the variables and the logical symbols, the formulas use the following symbols.
\begin{itemize}
\item In the valued field sort:
the language of rings.
\item In the residue field sort: the language of rings.
\item In the $\Z$-sort: the Presburger language.
\item a symbol $\ord(\cdot)$ for the valuation map from the nonzero elements of the valued field sort to the $\Z$-sort, and {a symbol} $\ac(\cdot)$ for the so-called angular component, which is a {multiplicative} function from the valued field sort to the residue field sort (more about this function below).
\end{itemize}

A formula in this language can be interpreted in any discretely valued field $F$ which comes with a
uniformizing element $\varpi$,
by letting the variables range over $F$, {over its} residue field $k_F$, and {over} $\Z$,
respectively, depending on the sort they belong to;
$\ord$ is the valuation map (defined on $F^\times$ and such that $\ord(\varpi)=1$), and $\ac$ is defined as follows.
If $x$ is a unit in $\cO_F$ (that is, $\ord(x)=0$), then $\ac(x)$ is the residue of $x$ modulo $\varpi$ (thus, an element of the residue field). For all other
nonzero $x$, one puts $\ac(x) :=  \varpi^{-\ord(x)}x \bmod (\varpi)$; thus, $\ac(x)$ is the first non-zero coefficient of the $\varpi$-adic expansion of $x$. Finally, we define $\ac(0)=0$.

Thus, a formula $\varphi$ in this language with $n$ free valued-field variables, $m$ free residue-field variables, and $r$ free $\Z$-variables gives naturally, for each discretely valued field $F$, a subset $\varphi(F)$ of
$F^n\times \rf_F^m\times \Z^r$: namely, $\varphi(F)$ is the set of all the tuples where the interpretation of $\varphi$ in $F$ is ``true'' (where a variable appearing inside a quantifier runs over either $F$, $k_F$, or $\ZZ$, respectively, depending on its sort).

\subsubsection{Adding constants to the language}
\label{sect:ldpo}

In many situations one needs to work with geometric objects defined over some fixed base number field, or over its ring of integers $\ri$. In such situations,
on top of the symbols for the constants that are already present in the above language (like $0$ and $1$), we will add to the Denef-Pas language all elements of
$\ri[[t]]$ as extra symbols for constants in the valued
field sort.
We denote the resulting language by $\ldpo$.

Given a complete discretely valued field $F$    {which} is an algebra over $\ri$ via {a chosen}
algebra homomorphism $\iota:\ri\to F$,    {with} a chosen uniformizer $\varpi$ of the valuation ring $\cO_F$ of $F$, one
can interpret the formulas in $\ldpo$ as described in the previous subsection, where
the new constants from $\ri$ are interpreted as elements of $F$ by using $\iota$, the constant symbol $t$ is interpreted as the uniformizer $\varpi$, and thus, by the completeness of $F$, elements of $\ri[[t]]$ can be naturally interpreted in $F$ as well.

\begin{defn}\label{AO}
Let $\cC_\ri$ be the collection of all triples $(F, \iota , \varpi)$, where $F$ is a non-Archimedean local
field which allows at least one ring homomorphism from $\ri$ to $F$,  the map $\iota:\ri\to F$ is such a ring homomorphism, and $\varpi$ is a uniformizer for $F$.
Let $\cA_\ri$ be the collection of triples $(F, \iota , \varpi)$ in $\cC_\ri$
with
$F$ of characteristic zero, and let $\cB_\ri$ be the collection of the triples $(F, {\iota} , \varpi)$ such that $F$ has positive characteristic.

Given an integer $M$, let $\cC_{\ri, M}$ be the collection of $(F,{\iota},\varpi)$ in $\cC_\ri$ such that the residue field of $F$ has characteristic larger than $M$, and similarly for $\cA_{\ri, M}$
and $\cB_{\ri, M}$.
\end{defn}

Since our results and proofs are independent of the choices of the map $\iota$ and the
uniformizer $\varpi$, we will often just write $F\in \cC_\ri$, instead of naming the whole triple.
For any $F\in \cC_{\ri}$, write $\cO_F$ for the valuation ring of $F$, $k_F$ for its  residue field, and $q_F$ for the cardinality of $k_F$.

\subsection{Definable sets and motivic  functions}\label{subsub:functions}

{The $\ldpo$-formulas introduced in the previous section allow us to obtain a} field-independent notion of subsets of $F^n\times k_F^m\times \ZZ^r$ for $F\in \cC_\ri$.

\begin{defn}\label{defset}
A collection
$X = (X_F)_{F}$ of subsets $X_F\subset F^n\times \rf_F^m\times \Z^r$ for $F$ in $\cC_\ri$ is called a \emph{a definable set} if there is an
$\ldpo$-formula $\varphi$ such that $X_F = \varphi(F)$ for each $F$, where $\varphi(F)$ is as explained at the end of \S \ref{DP}.
\end{defn}

{By Definition \ref{defset}}, a ``definable set'' is actually a collection of sets indexed by $F\in \cC_\ri$; such practice is not uncommon in model theory and has its analogues in algebraic geometry.
{A particularly simple definable set is $(F^n\times k_F^m\times \ZZ^r)_F$; we introduce
the simplified notation}
${\rm VF}^n\times {\rm RF}^m\times \ZZ^r$, where ${\rm{VF}}$ stands for valued field and ${\rm{RF}}$ for residue field.
{We apply the typical set-theoretical notation to definable sets $X, Y$,    {e.g.},
$X \subset Y$ (if $X_F \subset Y_F$ for each $F \in \cC_\ri$), $X \times Y$, and so on.}

For definable sets $X$ and $Y$, a collection $f = (f_F)_F$ of functions $f_F:X_F\to Y_F$ for $F\in\cC_\ri$ is called a definable function and denoted by $f:X\to Y$ if
the collection of graphs of the $f_F$ is a definable set.

There are precise quantifier elimination statements for {the} Denef-Pas language,
which  lead to the study of the geometry of definable sets and functions for $F$ in $\cC_{\ri,M}$, where often $M$ is large and depends on the data defining the set. To give an overview of these results is beyond of the scope of  this survey, but the reader may consult \cite{CCL:metric}, \cite{CCL:lipschitz}, \cite{CH: Lipschitz} and others.

We now come to motivic  functions, for which definable functions are the building blocks as mentioned above.
\begin{defn}\label{motfun}
Let $X = (X_F)_F$ be a definable set.
A collection $f = (f_F)_F$ of functions $f_F:X_F\to\CC$ is called \emph{a motivic function} on $X$ if and only if
there exist integers
$N$, $N'$, and $N''$, such that, for all $F\in \cC_\ri$, 
$$
f_F(x)=\sum_{i=1}^N   q_F^{\alpha_{iF}(x)} ( \# (Y_{iF})_x  )  \big( \prod_{j=1}^{N'} \beta_{ijF}(x) \big) \big( \prod_{\ell=1}^{N''} \frac{1}{1-q_F^{a_{i\ell}}} \big), \mbox{ for } x\in X_F,
$$
for some
\begin{itemize}
\item nonzero integers $a_{i\ell}$, 
\item definable functions $\alpha_{i}:X\to \ZZ$ and $\beta_{ij}:X\to \ZZ$, \item definable sets  $Y_i\subset X\times {\rm RF}^{r_i}$,
\end{itemize}
where, for $x\in X_F$,  $(Y_{iF})_x$ is the finite set $\{y\in k_F^{r_i}\mid (x,y)\in Y_{iF}\}$.
\end{defn}

\subsection{Motivic exponential functions}\label{subsub:expfunctions}

For any local field $F$, let $\scD_F$ be the set of the additive characters on $F$ that are trivial on    {the maximal ideal $\mathfrak m_F$} and nontrivial on $\cO_F$.
{For $\Lambda \in \scD_F$,
we will write $\bar\Lambda$ for the induced additive character on $k_F$.

\begin{defn}\label{expfun}
Let $X  = (X_F)_F$ be a definable set.
A collection $f = (f_{F,\Lambda})_{F,\Lambda}$ of functions $f_{F,\Lambda}:X_F\to\CC$ for $F\in \cC_\ri$ and $\Lambda\in \scD_F$ in is called \emph{a motivic exponential function} on $X$ if 
there exist integers $N>0$ and $r_i\geq 0$, motivic functions $f_i=(f_{iF})$, definable sets $Y_i\subset X\times {\rm RF}^{r_i}$ and definable functions $g_i:Y_i\to {\rm VF}$ and $e_i:{Y_i}\to {\rm RF}$ for $i=1,\ldots,N$, such that for all $F\in \cC_\ri$ and all $\Lambda\in \scD_\K$ 
$$
f_{F,\Lambda}(x)=\sum_{i=1}^N   f_{iF}(x)\Big( \sum_{y \in (Y_{iF})_x}\Lambda\big(g_i(x,y)\big) \cdot \bar \Lambda\big(e_i(x,y)\big)\Big)\mbox{ for all } x\in X_F.
$$
\end{defn}

{Compared to Definition \ref{motfun}, the counting operation $\#$ has been replaced by taking exponential sums, which makes the motivic exponential functions a richer class than the motivic functions.}
Motivic (exponential) functions form classes of functions which are stable under integrating some of the variables out, see Theorem \ref{thm:mot.int.} below. In fact, the search for this very property led to their definition.

\subsection{Families}
   {There are natural notions of ``uniform families'' of the objects introduced in Definitions~\ref{defset},
\ref{motfun} and \ref{expfun}. We fix the following terminology.}

\begin{defn}\label{deffam}
Suppose that $A$ is a definable set.
A {\emph{family of definable sets}} with parameter in $A$ is a definable set $X\subset A\times Y$ with some definable set $Y$.
We denote such a family by $(X_a)_{a \in A}$.
For each $F\in \cC_\ri$ and $a\in A_F$, the sets $X_{F,a}:= \{y\in Y_F\mid (a,y)\in X_F  \}$
are called the family members of the family $X_F = (X_{F,a})_{a\in A_F}$.

Similarly, for $(X_a)_{a \in A}$ as above, an (exponential) motivic, resp.~definable, function $f$ on $X$
can be considered as a \emph{{family $(f_a)_{a \in A}$} of (exponential) motivic}, resp.~\emph{definable
functions} with parameter in $A$. Its members are functions $f_{F,a}$ or $f_{F,\Lambda,a}$ for $F\in \cC_\ri$, $\Lambda\in \scD_F$ and $a\in A_F$.
\end{defn}

Whenever we call a specific function $f_{F}$ motivic, we implicitly think of $F$ as varying 
and we mean that there exists an $M > 0$ and a
motivic function    {$g$ with $f_F = g_F$} for all $F \in C_M$.
Likewise, we think of $F$ as varying whenever a specific set $X_{F}$ or function $f_{F}$ 
is called definable.
All the results presented in this article will only be valid for sufficiently big residue characteristic.
In particular, for any of the uniform objects $X$ we introduced in Defitions~\ref{defset}, \ref{motfun}, \ref{expfun},
we are actually only interested in $(X_F)_{F \in \cC_{\ri,M}}$ for $M$ sufficiently big.}

\subsection{Measure and integration}

To integrate a motivic function $f$ on a definable set $X$, we need a uniformly given family of measures on each $X_F$. For $X = {\rm VF}$, we put the Haar measure on $X_F = F$ so that $\cO_F$ has measure $1$; on
$k_F$ and on $\ZZ$, we use the counting measure. To obtain measures on
arbitrary definable sets $X$, we use ``definable volume forms'', as follows.

In a nutshell, a definable volume form $\omega$ on $X$ is a volume form defined by means of definable functions on definable coordinate charts, see \cite{cluckers-loeser-nicaise}.
Such a definable volume form specializes to a volume form $\omega_F$ on $X_F$ for every $F$
   {with sufficiently large residue characteristic; more precisely, a definable volume form specializes to a top degree differential form for every $F$, but its non-vanishing can be assured only when the residue characteristic of $F$ is sufficiently large.}

Any $F$-analytic subvariety of $F^n$, say, everywhere of equal dimension, together with an $F$-analytic volume form, carries a natural measure associated to the volume form, cf.~\cite{Bour}. This carries over to the
motivic setting, i.e., a volume form $\omega_F$ as above induces a measure on $X_F$. By a
``motivic measure'', we mean a family of measures on each $X_F$ that arises in this way from a definable
volume form.

Any algebraic volume form defined over $\ri$ on a variety $V$ over $\ri$ gives rise to a definable
volume form on $V(F)$ for    {$F\in \cC_{\cO, M}$ for a sufficiently large constant $M$}, in the above sense.
However, often one deals with volume forms on, say, orbits of elements of a group defined over a local field, or other volume forms defined over a local field, for which it is not always easy to obtain that they can be defined in a way which does not depend on the particular local field.

The classes of motivic (exponential) functions are natural to work with for the purposes of integration, as shown by the following theorem, which generalizes a similar result of
\cite{cluckers-loeser:fourier} to the class of functions summable in the classical sense of $L^1$-integrability.

\begin{thm}\cite{CGH-1}*{Theorem 4.3.1}\label{thm:mot.int.}
Let $f$ be a motivic  function, resp.~a motivic exponential function, on $X\times Y$ for some definable sets $X$ and $Y$, with $Y$ equipped with a motivic measure $\mu$. Then there exist a motivic  function, resp.~a motivic exponential function, $g$ on $X$ and an integer $M>0$ such that for each $F\in \cC_{\ri,M}$, each $\Lambda\in \scD_F$
and for each $x\in X_F$ one has
$$
g_F(x) = \int_{y\in Y_F}f_F(x,y)\,d\mu_F, \mbox{ resp. } g_{F,\Lambda}(x) = \int_{y\in Y_F}f_{F,\Lambda}(x,y)\,d\mu_F,
$$
whenever the function $Y_F\to\CC:y\mapsto f_F(x,y)$, resp.~$y\mapsto f_{F,\Lambda}(x,y)$ lies in $L^1(Y_F, \mu_F)$.
\end{thm}

\section{New transfer principles for motivic exponential functions}\label{sec:trpr}
Theorem \ref{thm:mot.int.} leads to the transfer of integrability and transfer of boundedness
principles, which we quote from \cite{CGH-1}.
(For simplicity, we quote the version without parameters, that is,
we take the parameter space to be a point in \cite{CGH-1}*{Theorem 4.4.1} and
\cite{CGH-1}*{Theorem 4.4.2}).
\begin{thm}\cite{CGH-1}*{Theorem 4.4.1}\label{thm:ti}
Let $f$ be a motivic  exponential
function on ${\rm VF}^n$.
Then there exists $M>0$, such that for the
fields $F\in\cC_{\ri, M}$,
the truth of the statement that $f_{F,\Lambda}$
is (locally) integrable
for all $\Lambda \in \scD_F$
depends only on the isomorphism class of the residue field of $F$.
\end{thm}

\begin{thm}\cite{CGH-1}*{Theorem 4.4.2}\label{thm:bounded}
Let $f$ be a motivic exponential
function on ${\rm VF}^n$.
Then, for some $M>0$, for the fields $F\in\cC_{\ri, M}$,
the truth of the statement that $f_{F, \Lambda}$ is (locally) bounded
for all $\Lambda \in \scD_F$
depends only on the isomorphism class of the residue field of $F$.
\end{thm}

In the versions with parameters of these results,    {$f$} is a family of motivic exponential functions
with parameter in a definable set $A$. In that case, each of the two theorems states that there exists a
motivic exponential function $g$ on $A$ such that the    {statement under consideration} is true for
$f_{F,\Lambda,a}$ iff $g_{F,\Lambda}(a) = 0$, that is, we relate the \emph{locus of integrability}, resp. the \emph{locus of boundedness} of the motivic (exponential) function $f$ to the zero locus of another motivic (exponential) function $g$.

\section{New ``uniform boundedness'' results for motivic functions}\label{sec:ub}

It turns out that the field-independent description of motivic functions leads not only to
transfer principles, but also to uniform (in the residue characteristic of the field) bounds for those motivic functions whose specializations are known to be bounded for individual local fields (for now, for motivic functions without the exponential).  We quote the following theorems from \cite{S-T}.

\begin{thm}\label{thm:presburger-fam}\cite{S-T}*{Theorem B.6}
Let $f$ be a motivic  function on $W \times \Z^n$,
where $W$ is a definable set. Then
there exist integers $a,b$ and $M$ such that for all $F\in \cC_{\ri, M}$ the following holds.

If there exists any (set-theoretical) function $\alpha:\ZZ^n\to\R$ such that
$$
| f_F(w,\lambda) |_{\R} \le \alpha(\lambda) \mbox{ on } W_F \times \Z^n,
$$
then one actually has
$$
| f_F(w,\lambda) |_{\R} \le q_F^{a+b \|\lambda\|} \mbox{ on } W_F \times \Z^n,
$$
where    {$\|\lambda\| = \sum_{i=1}^n |\lambda_i|$}, and where $|\cdot|_\R$ is the norm on $\R$.
\end{thm}

We observe that in the case with $n=0$,
the theorem yields that if a motivic  function $f$ on $W$ is such that $f_F$ is bounded on $W_F$ for each $F\in \cC_{\ri, M}$, then the bound for $|f_F|_\R$  can be taken to be $q_F^a$ for some $a\geq 0$, uniformly in $F$ with large residue characteristic.

The following statement allows one to transfer bounds which are known {for local fields of}
characteristic zero to local fields of positive characteristic, and vice versa.

\begin{thm}\label{thm:transfer-fam}\cite{S-T}*{Theorem B.7}
Let $f$ be a  motivic function on $W \times \Z^n$,
where $W$ is a definable set, and let $a$ and $b$ be integers. Then there exists $M$ such that, for any $F\in \cC_{\ri, M}$,
whether the statement
\begin{equation}\label{transfer}
f_F(w,\lambda)\le q_F^{a+b \|\lambda\|} \mbox{ for all } (w,\lambda) \in W_F \times \Z^n
\end{equation}
holds or not, only depends on the isomorphism class of the residue field of $F$.
\end{thm}

\subsection{A few words about proofs.}

The common strategy in the proofs of Theorems~   {\ref{thm:mot.int.}} -- \ref{thm:transfer-fam} (assuming the set-up of motivic integration developed in \cite{cluckers-loeser} and 
\cite{cluckers-loeser:fourier}) is the following.
The first step is to reduce the problem from
arbitrary motivic exponential functions $f$ to motivic exponential functions $g_i$ of a simpler form,
by writing $f$ as a sum $\sum_i g_i$. For some of the results
(e.g. to obtain a precise description of the locus of integrability or boundedness of $f$),
one needs to prove that not too much cancellation can occur in this sum; this is in fact
one of the main difficulties in the proofs of Theorems~\ref{thm:ti} and \ref{thm:bounded}.

In the next step, we eliminate all the valued-field variables, possibly at the cost of introducing more residue-field and $\ZZ$-valued variables, which uses the Cell Decomposition Theorem.
Once we have a motivic function that depends only on the residue-field and value-group variables,
we again break it into a sum of simpler terms, and then get rid of the residue-field variables.
Finally, the question is reduced to the study of functions $h : \ZZ^r \to \R$, which are,
up to a finite partition of $\ZZ^r$ into sets    {which are} definable in
the Presburger language (see Section~\ref{Pres}),
sums of products of linear functions in $x \in \ZZ^r$ and of powers of $q_F$, where the power also depends linearly on $x$. For most $x$, the behaviour of such a sum $h$ is governed by a dominant term,
and for a single term, each of the theorems is easy to prove.
To control the locus where no single term of $h$ is dominant, we use the powerful result by Wilkie about
``o-minimality of $\R_{\mathrm{exp}}$''.

\section{Which objects arising in harmonic analysis are definable?}
In the remainder of the article,
when we say ``definable'', we mean, definable in
the Denef-Pas language $   {\ldpo[\Z]}$
(introduced in Section~\ref{sect:ldpo}). 
If one wishes to think specifically of reductive grouse defined over a given global field $\bF$ with the ring of integers $\ri$ and its completions, then one can use the language $\ldpo$
instead.
We write $\cA_M$, $\cB_M$, and $\cC_M$ for $\cA_{\Z, M}$, $\cB_{\Z, M}$ and $\cC_{\Z, M}$, respectively (see Definition~\ref{AO}).

\subsection{Connected reductive groups}
We recall that a split connected reductive group over a local field is determined by its
root datum $\Psi=(X_\ast, \Phi, X^\ast, \Phi^\vee)$. By an \emph{absolute} root datum of $\bG$ we will mean
the root datum of $\bG$ over the separable closure of the given local field; it is a quadruple of this form.

By an unramified root datum we mean a root datum of an unramified reductive group over a local field $F$, i.e. a quintuple $\xi=(X_\ast, \Phi, X^\ast, \Phi^\vee, \theta)$, where $\theta$
is the action of the Frobenius element of $F^\ur/F$ on the first four components of $\xi$.

Given an unramified root datum $\xi$,
it is shown in \cite{cluckers-hales-loeser}*{\S 4} that a connected reductive group
$\bG(F)$ with root datum $\xi$ is ``definable using parameters'', for all local 
fields $F\in \cC_M$, where $M$ is determined by the root datum. This result is extended to all connected reductive groups (not necessarily unramified) in \cite{CGH-2}*{\S 3.1}, \cite{S-T}*{B.4.3}.
More precisely, by ``definable using parameters'', we mean that $\bG(F)$ is a member of a family
(in the sense of Definition~\ref{deffam}) of definable connected reductive groups
with parameter in a suitable definable set $Z$. This parameter encodes information about the extension over which $\bG$ splits (see \cite{CGH-2}*{\S 3.1} for details). In the remainder of this article, every definable or motivic object in $\bG(F)$ is parametrized by $Z$, and all results are uniform in the family.

\subsection{Haar measure and Moy-Prasad filtrations} If $\bG$ is unramified, the canonical measure on $\bG(F)$ defined by Gross \cite{gross:motive} is motivic for $F\in \cC_M$ (with $M$ depending only on the root datum that determines $\bG$), see \cite{cluckers-hales-loeser}*{\S7.1}.
For a root datum defining a ramified group $\bG$,    {to the best of our knowledge, this is presently not known}.    {However, in \cite{S-T}*{Appendix B} we prove} that
there is a constant $c$ that depends only on the root datum, such that Gross' measure can be renormalized by a factor between $q^{-c}$ and $q^c$ {(   {possibly} depending on $F$)} to give a motivic measure on $\bG(F)$,
where $q=p^r$ is the cardinality of the residue field of $F$, and
$p$ is assumed large enough so that the group is tamely ramified.

In fact, this question is closely related to whether Moy-Prasad filtration subgroups
$\bG(F)_{x,r}$ are definable (here $x$ is a point in the building of $\bG(F)$, and $r\ge 0$). This question has been answered in a number of important special cases. Namely,
if $\bG$ is unramified over $F$, and $x$ is a hyper-special point, then
we prove in \cite{CGH-2}*{\S 3} that $\bG_{x, r}$ is definable for all $r\ge 0$.
The same proof (which relies on the split case and taking Galois-fixed points) works in the case when $\bG$ is ramified over $F$, $x$ is a special point,
and $r>0$, but fails for $r=0$ since Galois descent no longer works. This is exactly the reason that the question whether the canonical measure is motivic is open in the ramified case.
We also prove in \cite{CGH-2}*{\S3} that for $x$ an optimal point in the alcove (in the sense of the definition by Adler and DeBacker, \cite{adler-debacker:bt-lie}),
$\bG(F)_{x,r}$ is definable for $r>0$, and for all $r\ge 0$ if $\bG$ is unramified over
$F$, and therefore the $G$-domain
$\bG(F)_r$ is definable. We also prove the corresponding results for the Lie algebra.
   {We note that the first results on definability of $\bG(F)_{x,0}$ in some special cases appeared in \cite{diwadkar:thesis}.}

\subsection{Orbital integrals and their Fourier transforms}
There are two ways to think of an orbital integral: as a distribution on the space $C_c^\infty(\bG(F))$ of locally constant compactly supported functions on $\bG(F)$ (or, respectively, on the Lie algebra), or to fix a test function $f$ on $\bG(F)$ (respectively, on $\fg(F)$), and consider the orbital integral
$\cO_\gamma(f)$ (respectively, $\cO_X(f)$) as a function on $\bG(F)$ (respectively, on the Lie algebra $\fg(F)$).

The earliest results concerning definability of orbital integrals
use the second approach:
in \cite{cunningham-hales:good}, C. Cunningham and T.C. Hales prove that
with $f$ fixed to be the characteristic function of $\fg(\ri_F)$
and {if} $\bG$ {is} 
unramified,
$\cO_X$ is a motivic function on the set of so-called \emph{good} regular elements in $\fg(F)$. In fact, they prove a slightly more refined result, giving explicitly the residue field parameters that control the behaviour of this function.
In \cite{hales:orbital_motivic}, T.C. Hales proves that orbital integrals are motivic ``on average''.
Both of these papers pre-date the most general notion of a motivic function, and because of this, extra complications arise from having to explicitly track the parameters needed to define orbital integrals.

{In order to be able to use} the first approach to orbital integrals (as distributions), we 
   {define a notion of a ``motivic distribution'' (without this name, it first appears in this context in \cite{cluckers-cunningham-gordon-spice}*{\S 3})}. We say that a distribution $\Phi$ is motivic if    {it is, in fact, a collection of distributions defined for each $F\in \cC_{M'}$ for some $M'>0$, and} for any
   {suitable} family $(f_a)_{a\in S}$ of motivic test functions 
indexed by a parameter $a$ varying over a definable set $S$,
there exists a motivic function {$g$} on $S$ {and an    {$M \ge M'$}} such that 
$\Phi(f_{a,F})=g(a)$ for all $F \in \cC_M$ and all    {$a \in S_F$,}
   {where by ``suitable'' we mean that $f_{F,a}$ is in the domain of 
$\Phi$ for all $F\in \cC_M$ and  $a\in S_F$.}

The most general current results on orbital integrals are:
\begin{enumerate}
\item Cluckers, Hales and Loeser prove in \cite{cluckers-hales-loeser} that
for $\bG$ -- unramified over $F$, and with $f$ fixed to be the characteristic function of
$\fg(\ri_F)$,
$X\mapsto \cO_X(f)$ is a motivic function on the set of regular semisimple elements in $\fg(F)$.
In fact, their method shows also that the orbital    {integrals of 
regular semisimple elements $X\in \fg(F)$ form a family of} 
motivic distribution in the sense defined above,    {indexed by $X$}.

\item Still under the assumption that $\bG$ is unramified, in \cite{S-T}*{\S B.6}, we prove that with a suitable normalization of measures,
non-regular orbital integrals 
{$\cO_\gamma$}
of tame elements (note that all elements are tame when the residue characteristic $p$ is large enough),  form a family of motivic distributions
(i.e., they  also depend on $\gamma\in \bG(F)$
in a definable way, thus showing that the orbital integrals are motivic in the sense of both approaches).

If the centralizer $C_G(\gamma)$ of $\gamma$
splits over an unramified extension, our normalization of the measure
coincides with the canonical measure on the orbit that one gets from the
canonical measure on $C_G(\gamma)$; if $C_G(\gamma)$ is ramified, the motivic and canonical  measures might differ by a constant not exceeding a fixed power of $q$.
This result should easily generalize to the case when $G$ is ramified as well (though one more constant from the measure on $G$ would appear).

We also note that even though we state the result for orbital integrals on the group, since this is what is required for the estimates that are the main goal of \cite{S-T}*{Appendix B}, in the course of the proof we essentially prove the analogous statement for the Lie algebra as well.

\item It follows from the results of McNinch \cite{mcninch:nilpotent} and Jantzen \cite{jantzen} that the set of nilpotent elements $\mathcal N$ in the Lie algebra $\fg(F)$ is definable,
cf \cite{CGH-2}*{\S 4.1}.
We prove  in
\cite{CGH-2}*{\S 4}
that nilpotent orbital integrals (on the Lie algebra)    {form a family of} motivic distributions,    {indexed by the set of nilpotent elements $\mathcal N$}.

\item On the Lie algebra, one can consider the Fourier transform
$\widehat{\cO}_X(f)=\cO_X(\widehat f)$ of an orbital integral. By a theorem of Harish-Chandra
in the characteristic zero case, and of Huntsinger \cite{adler-debacker:mk-theory} in positive characteristic, this distribution is represented by a locally constant function supported
on the set of regular elements $\fg(F)^{\reg}$. We will denote this function by
$\widehat{\mu}_X(Y)$.
In \cite{CGH-2}*{Theorem 5.7, Theorem 5.11} we prove that:
\begin{enumerate}
\item For nilpotent elements $X$, $\widehat \mu_X(Y)$ is a motivic
function defined on $\mathcal N\times \fg$.
\item For regular semisimple elements $X$,  $\widehat \mu_X(Y)$ is a motivic
function defined on $\mathcal \fg^\reg\times \fg$.
\end{enumerate}
\end{enumerate}

\subsection{Classification of nilpotent orbits,  and Shalika germs}
It is well-known that when the residue characteristic of the field $F$ is large enough,
there are finitely many nilpotent orbits in $\fg(K)$.
Using Waldspurger's parametrization of nilpotent orbits for classical Lie algebras
\cite{waldspurger:nilpotent}*{\S I.6}, it is possible to prove that these orbits are parametrized by the residue-field points of a definable set.
This approach is carried out in detail for odd orthogonal groups in \cite{diwadkar:o_n}.
Alternatively, one can use DeBacker's parametrization of nilpotent orbits that uses Bruhat-Tits theory. This approach is demonstrated for the Lie algebra of type  $G_2$ (as well as for $\fsl_n$) in J. Diwadkar's thesis \cite{diwadkar:thesis}.

Using both parametrizations of nilpotent orbits, and the explicit matching between them
due to M. Nevins, \cite{nevins:param}, L. Robson proves
in \cite{lance:thesis}
that for $\fsp_{2n}$, Shalika germs are motivic functions, up to a ``motivic constant''.
More precisely,    {one of the main results of \cite{gordon-robson} is}:
\begin{theorem} Let $\fg =\fsp_{2n}$. Let $\mathcal N$ be the set of nilpotent elements in $\fg$.  Then
\begin{enumerate}
\item There exists a definable set
$\mathcal E\subset {\rm RF}^n$ (that is, having only residue-field variables), and a definable  function $h:\mathcal N\to \mathcal E$, such that
for every $d\in \mathcal E$, $h^{-1}(d)$ is an adjoint orbit, and each orbit appears exactly once as the fibre of $h$.
\item There exist a finite collection of integers  ${(a_i)}_{i=0}^{r}$,  a motivic function $\Gamma$ on
$\mathcal E\times\fg^{\reg}$, and a constant $M>0$, such that for $F\in \cC_M$,
for every $d\in {\mathcal E}_F$, the function 
$\left(q^{a_0}\prod_{i=1}^r \frac 1{1-q^{a_i}}\right)^{   {-1}}\Gamma_F(d, \cdot)$ is a representative of the Shalika germ on $\fg^{\reg}$ corresponding to $d$.
\end{enumerate}
\end{theorem}

In a different context of groups over $\C\llp t\rrp $, a motivic interpretation of the
sub-regular germ is given in E. Lawes' thesis, \cite{lawes:thesis}.

\subsection{Transfer factors and other ingredients of the Fundamental Lemma}
In \cite{cluckers-hales-loeser}, the authors prove that transfer factors, weights, and other ingredients of the Fundamental Lemma are motivic functions.

\subsection{Characters}
So far, the only results on Harish-Chandra characters of representations have been restricted to the cases where the character can be related to an orbital integral.
One of the fundamental problems is that the parametrization of supercuspidal representations via
$K$-types (and the classification of the $K$-types due to J.-K. Yu and Ju-Lee Kim) uses multiplicative characters of the field, which cannot be easily incorporated into the motivic integration framework. Because of this, at present we do not have a way of talking about representations in a completely field-independent way.
There are two partial results on characters of some representations in a suitable neighbourhood of the identity, where Murnaghan-Kirillov theory holds.
\begin{enumerate}
\item The main result of \cite{gordon:depth0} can be expressed in modern terms as:
\emph{the character of a depth-zero representation of $\mathrm{Sp}_{2n}$ or
$\mathrm{SO}_{2n+1}$ is a motivic distribution (in the sense defined above) on the set of topologically unipotent elements}.
\item For positive-depth so-called \emph{toral, very supercuspidal} representations,
it is shown in  \cite{cluckers-cunningham-gordon-spice}, that up to a constant, their Harish-Chandra characters are motivic exponential functions on a neighbourhood of the identity on which Murnaghan-Kirillov theory works, and an explicit parameter space (over the residue field) is constructed for these representations up to crude equivalence, where we identify two representations if their characters on this neighbourhood of the identity are the same.
\end{enumerate}

However, thanks to Harish-Chandra's local character expansion, extended to an expansion
near an arbitrary tame semisimple element by J. Adler and J. Korman 
\cite{adler-korman:loc-char-exp},
one can prove results about characters, using motivic integration, even  without
knowing that they are themselves motivic functions. We give an example of such an application in the next section.

We note that proving that the coefficients of Harish-Chandra's local character expansion
are motivic is a question of the same order of difficulty as the question about the characters themselves, and so it is presently far from known.

\section{Local integrability results in positive characteristic}
Following Harish-Chandra (and using the modification by R. Kottwitz
\cite{S-T}*{Appendix A} for non-regular elements),
for a semisimple element $\gamma\in \bG(F)$, we define the discriminant
$$D^G(\gamma):=
\prod_{\alpha\in \Phi\atop{ \alpha(\gamma)\neq 1}}|1-\alpha(\gamma)|, \quad \text{and}
\quad D(\gamma): =\prod_{\alpha\in \Phi}|1-\alpha(\gamma)|,
$$
where $\Phi$ is the root system of $\bG$.

On the Lie algebra, define
$D(X)=\prod_{\alpha\in \Phi}|\alpha(X)|$.
Following \cite{kottwitz:clay}, we call a function $h$ defined on $\fg(F)^\reg$
\emph{nice}, if it satisfies the following conditions:
\begin{itemize}
\item when extended by zero to
all of $\fg(\K)$, it is locally integrable, and
\item the function $D(X)^{1/2}h(X)$ is locally bounded on $\fg(\K)$.
\end{itemize}

Similarly, call a function on $\bG(\K)^\reg$ \emph{nice}, if it satisfies
the same conditions on $\bG(\K)$, with $D(X)$ replaced by its group version
$D(\gamma)$.

Combining the theorems of Harish-Chandra in characteristic zero with the Transfer Principles of \S \ref{sec:trpr} above, in \cite{CGH-2} we  obtain
\begin{thm}\label{thm:orb int loc int}
\begin{enumerate}
\item
For a connected reductive group $\bG$ (defined via specifying a root datum), there exists a constant $M_\bG>0$ that depends only on the absolute root datum of $\bG$,
such that
for every $\K\in \cB_{M_\bG}$, and every
nilpotent orbit $\cO$ in $\fg(\K)$,
the function $\widehat \mu_{\cO}$  is a nice function on $\fg(\K)$.

\item
There exists a constant $M_\bG^\reg >0$ that depends only on the absolute
root datum of $\bG$,
such that for $\K\in \cB_{M_\bG^\reg }$, for every $X\in \fg(\K)^\reg$,
the function $\widehat \mu_X$ is a nice function on $\fg(\K)$.
\end{enumerate}
\end{thm}

Local character expansion,    {which is crucial for our next result}, in large positive characteristic
is proved by DeBacker \cite{debacker:homogeneity} near the identity,
and by Adler-Korman \cite{adler-korman:loc-char-exp} near a general tame semisimple element.
These results require an additional hypothethis on the existence of the so-called mock exponential map $\mexp$, discussed in \cite{adler-korman:loc-char-exp}, which we will not quote here; we note only that for classical groups one can take $\mexp$ to be the Cayley transform.

\begin{thm}\label{thm:char}
Given an absolute root datum $\Psi$, there exists a constant $M>0$, such that for
every $F\in \cB_{M}$, every connected reductive group $\bG$ defined over $F$ with the absolute root datum $\Psi$, and every  admissible
representation $\pi$ of ${\bG}(F)$,
the Harish-Chandra character $\theta_{\pi}$ is a nice function on
$\bG(F)$, provided that mock exponential map $\mexp$ exists for $\bG(F)$.
\end{thm}

Finally, Theorem \ref{thm:orb int loc int} also implies (thanks to a result of DeBacker) that Fourier transforms of general invariant distributions on $\fg(F)$ with support bounded modulo conjugation are represented by nice functions  in a neighbourhood of the origin. 

\section{Uniform bounds for orbital integrals}

In order to state the next result we need to introduce the notation for a family
of test functions, namely, the generators of the spherical Hecke algebra. Assume that $\bG$ is an unramified reductive group with root datum $\xi$.
Let  $K={\underline G}({\ri_F})$ be a hyperspecial maximal compact subgroup (where
${\underline G}$ is the smooth model associated with 
a hyperspecial point by Bruhat-Tits theory).

We have a Borel subgroup $B=TU$, and let $A$ be the maximal $F$-split torus in $T$.
For $\lambda\in X_\ast(A)$, let $\tau_\lambda^G$ be the characteristic function of
the double coset $K\lambda(\varpi)K$.

In \cite{S-T}*{Appendix B}, we prove, using Theorem~\ref{thm:presburger-fam},
\begin{thm}\label{thm:main2}
Consider an unramified root datum $\xi$. Then there exist constants $M>0$,  $a_\xi$ and $b_\xi$ that depend only on $\xi$,  such that for each non-Archimedean local field $F$
with residue characteristic at least $M$, the following holds. Let $\bG$ be a connected reductive algebraic group over
$F$ with the root datum $\xi$. {Let $A$ be a maximal $F$-split torus in $\bG$, and let $\tau_\lambda^G$ be as above.} Then
for all {$\lambda\in X_\ast(A)$} with $\|\lambda\|\le \kappa$,
$$
|\oi_\gamma(\tau_\lambda^G)|\le {q_F}^{a_\xi+b_\xi\kappa} D^G(\gamma)^{-1/2}$$
{for all semisimple elements $\gamma\in \bG(F)$.}
\end{thm}

It is sometimes useful to reformulate theorems of this type  to make an assertion about
a group defined over a global field. This is possible because given a connected reductive group
over a global field $\bF$, there are finitely many possibilities for the root data of 
   {$\bG_v=\bG\otimes_{\bF}\bF_v$}, as $v$ runs over the set of finite places of $\bF$, 
and    {$\bF_v$ denotes the completion of $\bF$ at $v$}.
Here we give an example of one such reformulation, which proved to be useful in \cite{S-T} 
({we also refer to \S 7 of loc.cit. for a more detailed explanation of the setting}).
\begin{thm}\label{thm:main}
Let $\bG$ be a connected reductive algebraic group over $\bF$, 
   {with a maximal split torus $A$}.
There exist constants $a_G$ and $b_G$ that depend only on the
global integral model of $\bG$ such that    {for
all but finitely many places $v$},
for all $\lambda\in X_\ast(A_v)$ with $\|\lambda\|\le \kappa$, 
$$
|\oi_\gamma(\tau_\lambda^G)|\le q_v^{a_G+b_G\kappa} D^G(\gamma)^{-1/2}$$
for all semisimple elements $\gamma\in \bG(\bF_v)$,
where $q_v$ is the cardinality of the residue field of $\bF_v$.
\end{thm}

\end{document}